\title{Exit time tails from pairwise decorrelation in hidden Markov chains, with applications to dynamical percolation}
\author{Alan Hammond \and  Elchanan Mossel \and G\'abor Pete}
\documentclass[11pt]{article}

\usepackage{amsmath}
\usepackage{amssymb}
\usepackage{amsthm}
\usepackage{amsfonts}
\usepackage{graphicx}
\usepackage{mathrsfs}

\input epsf.sty
\input labelfig.tex

\usepackage{color}
\newcommand{\margin}[2]{\textcolor{magenta}{*}\marginpar{\vskip #1 cm \textcolor{magenta} {\it #2 }}}

\newif\iffigures\figurestrue
\figuresfalse

\newif\ifhyper\IfFileExists{hyperref.sty}{\hypertrue}{\hyperfalse}

\ifhyper\usepackage{hyperref}
\def\hitem#1#2{\item[\hypertarget{#1}{#2}]\expandafter\gdef\csname LBL#1ITM\endcsname{#2}}
\def\iref#1{\hyperlink{#1}{\csname LBL#1ITM\endcsname}}
\else
\def\hitem#1#2{\item[{#2}]\expandafter\gdef\csname LBL#1ITM\endcsname{#2}}
\def\iref#1{{\csname LBL#1ITM\endcsname}}
\fi

\newif\ifdraft
\drafttrue 
\long\def\comment#1{}
\long\def\old#1{}

\numberwithin{equation}{section}
\numberwithin{figure}{section}

\newtheorem{theorem}{Theorem}
\numberwithin{theorem}{section}
\newtheorem{corollary}[theorem]{Corollary}
\newtheorem{lemma}[theorem]{Lemma}
\newtheorem{proposition}[theorem]{Proposition}
\newtheorem{conjecture}[theorem]{Conjecture}
\newtheorem{question}[theorem]{Question}

\theoremstyle{remark}\newtheorem{definition}[theorem]{Definition}
\theoremstyle{remark}

\let\qqed=\qed
\def\QED{\qqed\medskip}
\let\qed=\QED
\newcommand{\ignore}[1]{}

\newcommand{\R}{\mathbb{R}}

\newcommand{\C}{\mathbb{C}}
\newcommand{\Z}{\mathbb{Z}}
\newcommand{\N}{\mathbb{N}}

\def\Quad{\mathcal Q}

\def\SLEkk#1/{$\mathrm{SLE}(#1)$}
\def\SLEr#1/{$\mathrm{SLE(\kappa;#1)}$}
\def\SLEkr#1;#2/{$\mathrm{SLE(#1;#2)}$}
\def\SLEk/{\SLEkk{\kappa}/}
\def\SLEtwo/{\SLEkk2/}
\def\SLE/{$\mathrm{SLE}$}
\def\SLEab/{\SLEkr 4; {a/\hco-1}, {b/\hco-1}/}
\def\Var{\mathrm{Var}}
\def\Ito/{It\^o}
\def \eps {\epsilon}
\def \P {{\bf P}}
\def\md{\mid}
\def\Bb#1#2{{\def\md{\bigm| }#1\bigl[#2\bigr]}}
\def\BB#1#2{{\def\md{\Bigm| }#1\Bigl[#2\Bigr]}}
\def\Bs#1#2{{\def\md{\mid}#1[#2]}}
\def\Pb{\Bb\P}
\def\Eb{\Bb\E}
\def\PB{\BB\P}
\def\EB{\BB\E}
\def\Ps{\Bs\P}
\def\Es{\Bs\E}

\def \p {{\partial}}
\def \E {{\bf E}}

\def \proof {{ \medbreak \noindent {\bf Proof.} }}
\def\proofof#1{{ \medbreak \noindent {\bf Proof of #1.} }}

\def\Cl{\mathscr{C}}

\def\Piv{\mathsf{Piv}}

\def\noopsort#1{}

\def\bl{\begin{lemma}}
\def\el{\end{lemma}}
\def\bth{\begin{theorem}}
\def\eth{\end{theorem}}
\def\bc{\begin{corollary}}
\def\ec{\end{corollary}}
\def\bcj{\begin{conjecture}}
\def\ecj{\end{conjecture}}
\def\bpr{\begin{proposition}}
\def\epr{\end{proposition}}
\def\bde{\begin{definition}}
\def\ede{\end{definition}}
\newcommand{\be}{\begin{eqnarray}}
\newcommand{\ee}{\end{eqnarray}}
\newcommand{\bes}{\begin{eqnarray*}}
\newcommand{\ees}{\end{eqnarray*}}

\usepackage{mathrsfs}
\def\One{\1}
\def\1{1\!\! 1}
\def\F{\mathcal F}
\def\CC{\mathcal C}

\def\AA{\mathcal{A}}
\def\EE{\mathcal E}

\newcommand{\supp}{{\mathsf{supp}}}
\newcommand{\coupling}{\mathbf{Q}}

\newcommand{\jumpch}{Y}
\newcommand{\pie}{{\empty}}
\def\llra{\longleftrightarrow}

\begin{document}
\maketitle

\begin{abstract}
Consider a Markov process $\omega_t$ at stationarity and some event $\CC$ (a subset of the state-space of the process). A natural measure of correlations in the process is the pairwise correlation $\Pb{\omega_0,\omega_t \in \mathcal{C}} - \Pb{\omega_0 \in \CC}^2$. A second natural measure is the probability of the continual occurrence event $\big\{ \omega_s \in \CC, \, \forall \, s \in [0,t] \big\}$.  
We show that for reversible Markov chains, and any event $\CC$, pairwise decorrelation of the event $\CC$ implies a decay of the probability of the continual occurrence event $\big\{ \omega_s \in \CC \, \forall \, s \in [0,t] \big\}$ as $t \to \infty$. We provide examples showing that our results are often sharp. 

Our main applications are to dynamical critical percolation. Let $\CC$ be the left-right crossing event of a large box, and let us scale time so that the expected number of changes to $\CC$ is order 1 in unit time. We show that the continual connection event has superpolynomial decay. Furthermore, on the infinite lattice without any time scaling, the first exceptional time with an infinite cluster appears with an exponential tail.
\end{abstract}

\section{Introduction}


We study the relationship between pairwise decorrelation of a specific event  and the decay rate of the probability of continual occurrence of the event in {\bf reversible Markov processes}.
 In particular, Theorem~\ref{t.main} below states that any decay of the pairwise correlations
$$
\Ps{\omega_0,\omega_t \in A}-\Ps{\omega_0\in A}^2
$$ for the process $(\omega_t)_{t\geq 0}$ in stationarity implies a comparable decay of the joint probability 
$$
\Ps{\omega_s \in A\text{ for all }0\leq s \leq t}\,.
$$ 

Given a Markov process $\omega_t$ on $S$ with stationary probability measure $\pi$, its time 1 Markov operator $(T_1f)(\omega):=\Es{f(\omega_1)\md\omega_0=\omega}$ on $f\in L^2(S,\pi)$ is a normal operator, hence it has a spectral decomposition, with $\mathsf{Spec}(T_1)\subseteq \{z:|z|\leq 1\} \subset \C$ and $z=1$ being an obvious eigenvalue. Its {\bf spectral gap} is  defined as $g:=\inf\big\{|1-\lambda|: \lambda \in \mathsf{Spec}(T_1)\setminus\{1\}\big\}$, while its {\bf absolute spectral gap} is $g_*:=1-\sup\big\{|\lambda|: \lambda \in \mathsf{Spec}(T_1)\setminus\{1\}\big\}$. 
It is well-known (and not hard to see) that $g_*>0$ is equivalent to having an exponential decay of correlations for {\bf any} function $f: S\longrightarrow \R$ with $\E_\pie[f]=0$: 
$$
\Es{f(\omega_0)f(\omega_t)}\leq (1-g_*)^t \, \E_\pie[f^2]\,,
$$ 
for the process at stationarity. We write $\P$ both for both the law $\pi$ on static configurations and for the measure of the process 
at stationarity. Similarly, $\E$ denotes expectation of a function with respect to one or other of these laws, depending on whether the function is defined on static or dynamic configurations
It is quite classical for this case (in fact, $g>0$ suffices), see \cite{AKSz}, \cite{AFWZ}, \cite[Theorem 3.6]{expanders} and \cite[Theorem 9.2.7]{AlonSpencer}, that for {\bf any} set $A\subset S$ with stationary measure bounded away from 1, the exit-time tail $\Ps{\omega_s \in A\text{ for all }0\leq s \leq t}$ is exponentially small in $t$, with an exponent depending on the spectral gap and on $\pi(A)$. The strongest such bound is the one in \cite{AlonSpencer}, with a generalization in \cite[Theorem 5.4]{MODRSS}. 

Our Theorem~\ref{t.main} is a generalization of these results for the case when we have a pairwise correlation decay not for any function, but only for being in a given $A$ --- which might happen on a much faster time-scale than the mixing time of the entire chain. In other words, our generalization concerns the {\bf hidden Markov chain} $\1_{\{\omega_t \in A\}}$.
Furthermore, our Theorem~\ref{t.AKSz} gives a generalization of \cite[Theorem 5.4]{MODRSS} in a different direction, by showing that, assuming a spectral gap, the exit-time tail from $A$ is exponentially decaying provided that the probability that $\omega_t$ is in $A$ at every moment of a fixed time interval is bounded away from one (which may be the case even if $\pi(A)$ is arbitrarily close to 1). 

The exponential exit-time tail for Markov chains with spectral gap 
(such as random walks on expander graphs) 
has many applications in computer science including derandomization of algorithms \cite[Section 3]{expanders} and noise sensitivity \cite{MODRSS}, 
suggesting that our results may prove useful from such points of view. Nevertheless, our initial motivation comes from the study of dynamical percolation on planar lattices, which is the natural time evolution of critical percolation in the plane, a central model of statistical mechanics; see \cite{HgPS, BKS, SS, GPS1, GPS2a, GPS2b, HPS} for the original papers, and \cite{St,GSt} for surveys. The implications of our results to dynamical percolation will be explained in Section~\ref{s.dynperc}. 

We now state our main results in detail.

\subsection{Exiting an event with some pairwise decorrelation}\label{ss.intro.exiting}

We will consider continuous or discrete time  Markov processes, $(\omega_t)_{t\in\R}$ or  $(\omega_t)_{t\in\Z}$, on some state space $S$, with some (not necessarily unique) stationary probability measure $\pi$; we will always consider the process run in stationarity, i.e., with $\omega_0\sim\pi$. For functions $f:S\longrightarrow \R$, consider the usual inner product $(f,g):=\E_\pie[fg]$, and the Markov operator $(T_tf)(\omega):=\Es{f(\omega_t) \md \omega_0=\omega}$.  Let $\CC$ be a static event (i.e., measurable with respect to~$\omega_0$), suppose that $\pi(\CC)=\Pb{\omega_0\in \CC}=p$, and let $f=\One_\CC$. The decay of correlations of $f$ in time is often quantified by the function $d:(0,\infty) \to [0,\infty)$ in one of the following two inequalities:
\begin{equation}
\Pb{\omega_0,\omega_t \in \CC}-\Ps{\omega_0\in\CC}^2 = (f,T_tf)-(\E_\pie f)^2 \leq
 d(t)\,\Var[f]\label{e.decorr1}
 \end{equation}
and
\begin{equation}
\Var[T_tf] = (T_tf,T_tf)-(\E f)^2 \leq d(2t)\,\Var[f]\,,
\label{e.decorr2}
 \end{equation}
for all $t\in[0,\infty)$. Of course, for reversible Markov processes, (\ref{e.decorr1}) is equivalent to~(\ref{e.decorr2}). We will consider the cases where the decay of  $d(t)$ as $t \to \infty$ is either polynomial or (stretched) exponential.
Sometimes, one has a sequence of Markov processes $(\omega^n_t)_{t\in\R}$, $n\in\N$, on larger and larger finite state spaces, with the time parameter coming from the original time of the process rescaled by a function of $n$. In this case, the bounds are understood uniformly in $n$.

\begin{theorem}\label{t.main} In the above setting, assuming (\ref{e.decorr2}), we have that
\begin{equation}\label{eqminform}
\PB{\omega_s \in \CC \ \forall s\in [0,t]} \leq \min_{k \in \N^+} \left\{ \left( \frac{p+1}{2} \right)^k + \frac{16 p}{(1-p)^2} d\left(\frac{2t}{k}\right)\right\}\,,
\end{equation}
and therefore 
$$
\PB{\omega_s \in \CC \ \forall s\in [0,t]} \leq \begin{cases} t^{-\alpha+o(1)} & \text{if } d(t)= \Theta(t^{-\alpha}),\\
\exp\big(-t^{\frac{\alpha}{1+\alpha}+o(1)}\big) &  \text{if } d(t)=\exp(-\Theta(t^{\alpha}))\,,\end{cases}
$$
as $t\to\infty$, where the $o(1)$ terms depend only on $p$, $\alpha$ and the constant factors implicit in the $\Theta(\cdot)$ notation.
\end{theorem}

Examples and questions of sharpness and of non-sharpness in Theorem~\ref{t.main} appear in Section~\ref{sec.ex}.
\medskip

\noindent{\bf Remark.} Note that we make no assumption of reversibility in Theorem~\ref{t.main}.
However, as we have noted, for a reversible Markov chain, we may replace the assumption (\ref{e.decorr2}) in the statement by (\ref{e.decorr1}), which is a more familiar form in which to express decorrelation in a Markov process. 
\medskip

\noindent{\bf Motivation.} As we mentioned above, our main motivation is dynamical critical percolation on planar lattices: site percolation on the triangular lattice or bond percolation on $\Z^2$. Let $\CC$ be the left-right crossing event of a large box, and let us scale time such that the expected number of changes to $\CC$ is order 1 in unit time. Theorem~\ref{t.main} implies that the continual connection event has superpolynomial decay. See Corollary~\ref{cor.LR}.

\ignore{
Firstly, Theorem~\ref{t.main} implies that for the {\bf scaling limit of dynamical percolation} the probability that the {\bf left-right crossing} in a square is maintained for time $t$ is superpolynomially small (i.e., is at most $C_K t^{-K}$, for any $K>0$). Our best lower bound is only $\exp(-c\,t)$, but still, this crossing event seems to be an enlightening example of how spectral arguments should really not yield an exponential upper bound.
} 

\subsection{Exiting events defined on time intervals, assuming a spectral gap}\label{ss.intro.AKSz}

We consider a continuous time Markov process semi-group $T_t = e^{tQ}$ on some state space $S$, reversible with respect to a probability measure $\pi$. Then the infinitesimal generator $Q$ is self-adjoint (reversible) and negative semi-definite with respect to the usual inner product given by $\pi$, and its spectrum is contained in $(-\infty,0]$. We will assume that $Q$ has a spectral gap $\delta>0$ around the obvious eigenvalue 0;  then $T_t$ has an absolute spectral gap $1-e^{-\delta t}$, and the process is ergodic.

Let $\Omega$ be the space of all paths $\omega: \R \longrightarrow S$ of the Markov process under the probability measure $\P$, and let $L^2(\Omega,\P)$ denote all $L^2$ integrable functions from $\Omega$ to $\R$. For a subset $I \subseteq R$ we denote by $\F_I$ the sigma algebra 
generated by $\{ \omega(t) : t \in I\}$.

\begin{theorem}\label{t.AKSz}
Suppose that the generator $Q$  is reversible with respect to a probability measure $\pi$ and has spectral gap $\delta > 0$. Let $k \in \N^+$, and let $a_i,b_i \in \R$, $0 \leq i \leq k$ satisfy $a_i < b_i$ for such $i$, as well as 
$b_i \leq a_{i+1}$ for $0 \leq i \leq k-1$; we also permit $a_0 = -\infty$
as well as $b_k = \infty$.

Let $\AA_0, \ldots, \AA_k$ be subsets of $\Omega$, each $\AA_i$ measurable with respect to
$\F_{[a_i,b_i]}$. Then
\[
\PB{\bigcap_{i=0}^{k} \AA_i} \leq \sqrt{\P[\AA_0]} \sqrt{\P[\AA_k]} \prod_{i=0}^{k-1} \left[ \sqrt{\P[\AA_i]} \sqrt{\P[\AA_{i+1}]}+e^{-\delta (a_{i+1} - b_i)} \left(1-\sqrt{\P[\AA_i]} \sqrt{\P[\AA_{i+1}]} \right) \right].
\]
In particular, suppose that $\AA$ is measurable with respect to $\F_{[a,b]}$ with $\P[\AA]=p$, and $[a_i,b_i]=[a+s_i,b+s_i]$. Setting $t_i = s_{i+1}-s_i - (b-a)$ for all $0\leq i <k$, then, provided that 
$t_i \geq 0$ for all such $i$, 
$$
\PB{\omega(t+s_i)_{t \in [a,b]} \in \AA \text{ for }i=0,1,\dots,k} \leq p \prod_{i=0}^{k-1} \big(p+e^{-\delta t_i}(1-p)\big).
$$
\end{theorem}

\noindent{\bf Motivation.} As we will show in Corollary~\ref{cor.radial}, this theorem implies for dynamical percolation on the infinite lattice that the probability that there is no {\bf exceptional time} in $[0,t]$ with the cluster of the origin being infinite is exponentially small. 
This corollary plays a significant role in \cite{HPS}.

\medskip

\noindent{\bf Acknowledgments.} We each benefitted greatly from conversations with Oded Schramm while working on this project and many others. We feel very fortunate to have known and to have worked with him. 
\medskip

We also thank Christophe Garban and Yuval Peres for useful conversations and remarks.

This work was started in the summer of 2007 at the Theory Group of Microsoft Research, Redmond. We are also grateful to the Fields Institute in Toronto, where some of the work by AH and GP was done in 2011. AH was supported principally by the EPSRC grant EP/I004378/1. EM was supported by NSF awards DMS 0548249 and DMS 1106999, DOD ONR grant N000141110140, ISF award 1300/08 and a Minerva award. GP was supported by an NSERC Discovery Grant at the University of Toronto, and an EU Marie Curie International Incoming Fellowship at the Technical University of Budapest.

\section{Proofs}\label{s.proofs}

\subsection{Exiting an event with some pairwise decorrelation}\label{ss.exiting}

We now prove Theorem~\ref{t.main}. 
\proof
Let $p < \lambda < 1$. Consider the static event
$$
A_s:=\Big\{\omega \in S: \Ps{\omega_s \in \CC \md \omega_0=\omega} < \lambda \Big\}\,.
$$
Note that for $s=0$ we have $\Pb{A_0^c} = p$, while for large $s$ one expects $\Pb{A_s^c}$ to be small.

Let $\tau=t/k$, for some $k\in \Z_+$. We claim that, for $m \geq 0$,
\begin{equation}\label{eqexplain}
 \Pb{\omega_0 \in A_\tau^c\cap \CC;\ \omega_{j\tau}\in A_\tau\cap \CC \text{ for } 1 \leq j \leq m}
 \leq  \lambda^{(m-1)\vee 0}\, \Pb{\omega_0 \in A_\tau^c}.
\end{equation}
We may prove this by induction on $m$, the cases where $m \in \{ 0,1 \}$ being trivial. For $m \geq 2$, writing 
$B_{m}$ for the event on the left-hand side of 
(\ref{eqexplain}),
we have that
 $\Pb{B_m} = \Pb{B_{m - 1}} q_m$, 
where $q_m$ is the conditional probability of 
$\omega_{m \tau} \in A_\tau \cap \CC$ given $B_{m - 1}$. 
Note that $q_m$ is at most 
$\Pb{\omega_{m \tau} \in \CC \md B_{m - 1}}$.
The conditional distribution of $\omega_{(m - 1)\tau}$ given 
$B_{m - 1}$ being supported on the event $\omega_{(m - 1) \tau} \in A_\tau$, it follows from the Markov property that 
$\Pb{\omega_{m \tau} \in \CC \md B_{m - 1}} \leq \lambda$.
Hence, the inductive hypothesis at $m -1$ implies this statement at $m$, giving~(\ref{eqexplain}).

By the same argument, we see that, for each $m \geq 0$,
\begin{equation}\label{eqexplainplus}
 \Pb{ \omega_{j\tau}\in A_\tau\cap \CC \text{ for } 0 \leq j \leq m}
 \leq  \lambda^m.
\end{equation}

We find then that
\begin{align}
\PB{\omega_s \in \CC \ \forall s\in [0,t]}
&\le \Pb{\omega_{j\tau} \in \CC \text{ for } 0\leq j\leq k}\nonumber\\
&\le \Pb{\omega_{j\tau} \in \CC \cap A_\tau \text{ for } 0\leq j\leq k} \; + \nonumber\\
& \qquad \sum_{\ell=0}^k \Pb{\omega_{\ell\tau}\in A_\tau^c\cap \CC;\ \omega_{j\tau}\in A_\tau\cap \CC \text{ for } \ell < j \leq k}\nonumber\\
&\leq \lambda^k + \sum_{\ell=0}^k \lambda^{(k-\ell-1)\vee 0}\, \Pb{\omega_{\ell\tau}\in A_\tau^c}\nonumber\\
&\leq \lambda^k + \frac{2 - \lambda}{1-\lambda}\, \Ps{A_\tau^c}\,.\label{e.twoterms}
\end{align}
In the third inequality, (\ref{eqexplainplus}) was used to bound the first term on its left-hand side, while the summand was bounded using (\ref{eqexplain}) with $m=k-\ell$ and stationarity.


We need to find now an upper bound on $\Pb{A_s^c}$ for $s$ large. 
By the definition of $A_s$,
$$
\Eb{\One_{A_s^c}\, T_s f} = \Eb{f(\omega_s) \md A_s^c} \, \Pb{A_s^c}
\geq \lambda\, \Pb{A_s^c}\,,
$$
where, as before, $f=\One_\CC$. On the other hand,
$$
\Eb{\One_{A_s^c}\, T_s f} = \Eb{\One_{A_s^c} p} + \EB{\One_{A_s^c}\, (T_sf-\E f)}\,.
$$
Putting these two things together,
$$
\EB{\One_{A_s^c}\,(T_sf-\E f)}
\geq (\lambda-p)\, \Pb{A_s^c}\,.
$$
Applying Cauchy-Schwarz to the left-hand side,
$$
\| \One_{A_s^c}\|_2 \,\|T_sf-\E f\|_2 \geq (\lambda-p)\, \Pb{A_s^c}\,,
$$
so that we obtain
$$
\Var[T_sf]^{1/2}=\|T_sf-\E f\|_2
\geq (\lambda-p)\, \Pb{A_s^c}^{1/2}\,.
$$
Therefore, (\ref{e.decorr2}) implies that
\be\label{e.Atauc}
\Pb{A_s^c} \leq \frac{p-p^2}{(\lambda-p)^2} \, d(2s).
\ee
Now take $s=\tau$, where recall that $\tau=t/k$ for some $k \in \Z_+$ that we will shortly specify. Plugging (\ref{e.Atauc}) into~(\ref{e.twoterms}), we find that
$$
\PB{\omega_s \in \CC \ \forall s\in [0,t]} \leq \lambda^{k} + \frac{p-p^2}{(\lambda-p)^2} \, \frac{2-\lambda}{1-\lambda} \, d(2t/k)\,.
$$
Setting $\lambda=(p+1)/2$ yields (\ref{eqminform}). 

For the case $d(t)=t^{-\alpha}$, setting
$k= \lfloor K\,\log t \rfloor$ for a suitable constant $K=K(\lambda,\alpha)$ makes both terms $t^{-\alpha}$, as desired.
For the case $d(t)=\exp(-t^{\alpha})$, we set $k= \lfloor t^{\beta} \rfloor$, and optimize the upper bound by letting $\alpha(1-\beta)=\beta$, i.e., choosing $\beta=\alpha/(1+\alpha)$, and we are done. 
\QED

\subsection{Exiting events defined on time intervals, assuming a spectral gap}\label{ss.AKSz}

In this subsection we prove Theorem~\ref{t.AKSz}. There are two main ideas. The first is that if a chain has a spectral gap, then the associated Markov operator will be a strict $L^2$-contraction not only on functions with zero mean, but also on any function whose support has a stationary measure bounded away from 1. The second idea is to use conditional expectation and the Markov property to extend the first idea to functions defined not on the state space, i.e., at individual times, but on time intervals. These two ideas are formalized in the following lemma.

\begin{lemma}\label{l.AKSz}
Suppose that the infinitesimal generator $Q$  is reversible with respect to a probability measure $\pi$ and has spectral gap $\delta > 0$. Let $\AA_1, \AA_2$ be subsets of $\Omega$ that are measurable with respect to 
$\F_{\leq a}$ and $\F_{\geq a+t}$, respectively. 
Let $P_1$ and $P_2$
be the corresponding projection operators on $L^2(\Omega,\P)$, i.e.,
$P_i f(\omega)=f(\omega)\1_{\AA_i}(\omega)$ for every function $f$ on $\Omega$.
Then
\[
\| P_1 P_2\| \leq
\sqrt{\P[\AA_1]} \sqrt{\P[\AA_2]}+e^{-\delta t} \left(1-\sqrt{\P[\AA_1]} \sqrt{\P[\AA_2]} \right),
\]
where the norm on the left is the operator norm for operators from
$L^2(\Omega, \P)$ into itself.
\end{lemma}

\proof
Since $P_1$ and $P_2$ are self-adjoint and commuting $P_1P_2$ is also self-adjoint. 
By the definition of norm, duality and the positivity of $P_i$,  
\begin{align*}
\| P_1 P_2\| &= \sup \big\{ (P_1P_2 f_2, f_1): \|f_i\|_2= 1,  f_i \geq 0 \big\}\\
&= \sup\big\{ (P_2 f_2, P_1 f_1) : \|f_i\|_2=1,  f_i \geq 0 \big\}.
\end{align*}
For $f_1,f_2$ such that $P_1 f_1 \neq 0, P_2 f_2 \neq 0$, define 
\[
f_1' = \frac{P_1 f_1}{\| P_1 f \|_2}, \quad 
 f_2' = \frac{P_2 f_2}{\| P_2 f_2 \|_2},
\]
so $\|f_i'\|_2 =1$ and $\supp(f_i) \subseteq \AA_i$, for $i=1,2$.  Since $\|P_if_i\| \leq1$ and since $P_1,P_2$ are idempotent, we have
 \[
 (P_2 f_2, P_1 f_1) \leq (f_2',f_1') =  (P_2 f_2', P_1 f_1'). 
 \] 
Therefore,
\[
\| P_1 P_2\| = \sup\Big\{ \E[f_1 f_2] :\|f_i\|_2=1,\, f_i \geq 0,\, 
\supp(f_i)\subseteq \AA_i\text{ for } i=1,2\Big\}.
\]

Given $f_1$ and $f_2$ as in the last equation, let $g_i : S \longrightarrow \R$ be defined by
 $g_1(x) := \Es{f_1(\omega) \md \omega(a) = x}$ and $g_2(x) := \Es{f_2(\omega) \md \omega(a+t) = x}$. 
 Clearly, $\E_\pie[g_i] = \E[f_i]$ and $\E_\pie[g_i^2] \leq \E[f_i^2] = 1$. 
 Since $f_1$ is $\F_{\leq a}$ measurable and $f_2$ is $\F_{\geq a+t}$ measurable, by the Markov property it follows that
\[
\E[f_1 f_2] = \Eb{ \Es{ f_1 f_2 \md \omega(a), \omega(a+t)}} = \Eb{g_1(\omega(a))\, g_2(\omega(a+t))} = \E_\pie[g_1 T_t g_2].
\]
Using the fact that the spectral gap of $T_t$ is $1-e^{-\delta t}$, it follows from Lemma~\ref{l.easyfact} below that the last expression is bounded by
\[
\E_\pie[g_1] \E_\pie[g_2] + e^{-\delta t} (1-\E_\pie[g_1] \E_\pie[g_2]) = 
\E[f_1] \E[f_2] + e^{-\delta t} (1-\E[f_1] \E[f_2]),
\]
and using $\E[f_i] = \E[f_i\, \1_{\AA_i}] \leq \|f_i\|_2 \|\1_{\AA_i}\|_2 \leq \sqrt{\P[\AA_i]}$ we obtain that
\[
\| P_1 P_2\| \leq 
\sqrt{\P[\AA_1]} \sqrt{\P[\AA_2]}+e^{-\delta t} \left(1-\sqrt{\P[\AA_1]} \sqrt{\P[\AA_2]} \right),
\]
as stated. 
\qed
\medskip

The proof of the previous lemma used the following easy fact.

\begin{lemma}\label{l.easyfact}
Let $M$ be an ergodic transition matrix for a
Markov chain on the set $S$ which is reversible with respect to
the probability measure $\pi$ and which has spectral gap $\delta>0$.
Let $g_1, g_2 : S \longrightarrow \R$ be two functions with $L^2$ norm at most one. 
Then 
\[
\E_\pie[g_1 M g_2] \leq \E_\pie[g_1] \E_\pie[g_2] + (1-\delta)\big(1-\big|\E_\pie[g_1] \E_\pie[g_2]\big|\big).
\]
\end{lemma}

\proof
Abbreviating $\One = \One_S$, we set
$h_i = g_i - \E_\pie[g_i] \1$. Then $h_i$ is orthogonal to the constant functions and $\1$ is a $1$-eigenvector of $M$. Therefore, 
\[
\E_\pie[g_1 M g_2] = \E_\pie[g_1] \E_\pie[g_2] + \E_\pie[h_1 M h_2].
\]
Using the spectral gap of $M$, we get:
\begin{align*}
\E_\pie[h_1 M h_2] \leq \| h_1 \|_2 \| M h_2 \|_2 \leq (1-\delta) \|h_1 \|_2 \| h_2 \|_2 
&\leq (1-\delta) \sqrt{(1-\E_\pie[g_1]^2)(1-\E_\pie[g_2]^2)} \\
&\leq (1-\delta) \big(1-\big|\E_\pie[g_1] \E_\pie[g_2]\big|\big), 
\end{align*}
where the last inequality follows from the inequality $(1-x^2)(1-y^2) \leq (1-xy)^2$, valid for all $x,y \in \R$, and we are done.
\qed

\proofof{Theorem~\ref{t.AKSz}}
 Let $P_i$ denote the projection onto $\AA_i$, as in
Lemma~\ref{l.AKSz}. 
It is easy to see that
\[
\PB{ \bigcap_{i=0}^k \AA_i} =
\EB{\1\, \big(\prod_{i=0}^k P_i\big)\, \1} =
\EB{\1\, \big(\prod_{i=0}^k P_i^2\big)\, \1},
\]
since the projection $P_i$ satisfies $P_i^2 = P_i$, and the order in which the projections act does not matter. These operators are also self-adjoint, so that 
\begin{align*}
\EB{\1 \, \big(\prod_{i=0}^k P_i^2\big)\, \1}=\EB{ (P_k\1) \,  P_k \, \big( \prod_{i=1}^{k-1} P_i^2 \big) \,  P_0 \, (P_0\1) }
&\leq \|\1_{\AA_0}\|_2\, \|\1_{\AA_k}\|_2 \, \Big\| P_k \, \big( \prod_{i=1}^{k-1} P_i^2 \big) \,P_0 \Big\| \\
& \leq \sqrt{\P[\AA_0]} \sqrt{\P[\AA_k]} \, \prod_{i=0}^{k-1} \| P_i P_{i+1} \|,
\end{align*}
where, in the rightmost expressions, the new notation denotes the operator norm from $L^2(\Omega,\P)$ to itself.
By Lemma~\ref{l.AKSz}, we have that, for each $i \in \{ 0, \dots, k-1 \}$,
\[
\|P_{i} P_{i+1}\| \leq \sqrt{\P[\AA_i]} 
\sqrt{\P[\AA_{i+1}]}+e^{-\delta (a_{i+1} - b_i)} \left(1-\sqrt{\P[\AA_i]} \sqrt{\P[\AA_{i+1}]} \right).
\]
Hence,
$$
 \prod_{i = 0}^{k-1} \| P_{i} P_{i+1} \| \leq
\prod_{i=0}^{k-1} \left[ \sqrt{\P[\AA_i]} \sqrt{\P[\AA_{i+1}]}+e^{-\delta  (a_{i+1} - b_i)} \left(1-\sqrt{\P[\AA_i]} \sqrt{\P[\AA_{i+1}]} \right) \right],
$$
and the proof is complete.
\qed

\section{Examples concerning Theorem \ref{t.main}}\label{sec.ex}

\subsection{An example where Theorem \ref{t.main} is sharp}

We give our example (of a reversible Markov chain) in terms of conductances on edges; see \cite[Chapter 2]{LPbook} for an exposition of the relevant theory.
Let $\beta \in  (1,2)$. Consider the graph with vertex set $\Z^* = \Z \setminus \{ 0 \}$
and edge set given by the nearest-neighbour edges among these vertices, plus the edge
$(-1,1)$, and a self-loop on each vertex. Equip the edges with conductances $c_{n,n+1} = c_{-(n+1),-n} := n^{-\beta}$ for each $n \geq 1$, $c_{-1,1} = c_{1,1}= c_{-1,-1}:=1/2$, and $c_{n,n}:=c_{n,n-1}+c_{n,n+1}$ for $|n|\geq 2$. Consider the discrete-time random walk $\jumpch$ on this graph equipped with this set of conductances. The sum of the conductances being finite, $\jumpch$ has a finite invariant measure $\pi$, whose value at a vertex is the sum of the conductances over incident edges; this is the unique invariant measure given its total mass.

Let $\P$ denote the law of $\jumpch$ run in stationarity.
Set $\CC = \N^+$, hence $\tfrac{\pi(\CC)}{\pi(\Z^*)}=1/2$. We will show that
\begin{equation}\label{eqctrwmin}
 \Pb{ \jumpch_s \in \CC \  \forall s \in [0,t] } = t^{\frac{1-\beta}{2} + o(1)}
\end{equation} 
and that
\begin{equation}\label{eqctrw}
 \Pb{ \jumpch_0, \jumpch_t \in \CC} -  \Pb{ \jumpch_0 \in \CC}^2 \leq \tfrac{1}{2} \Pb{ \jumpch_s \in \CC \  \forall s \in [0,t] }\,;
\end{equation}
moreover, for any $\epsilon > 0$, and for all $t > 0$ sufficiently high,
\begin{equation}\label{eqctrwplus} 
\Pb{ \jumpch_0, \jumpch_t \in \CC} -  \Pb{ \jumpch_0 \in \CC}^2 \geq \big( \tfrac{1}{2} - \epsilon \big) \Pb{ \jumpch_s \in \CC \  \forall s \in [0,t] }
\,.
\end{equation}
These show that Theorem~\ref{t.main} is basically sharp in the regime of polynomial decay; in fact, (\ref{eqctrwplus}) is a stronger bound than what follows from Theorem~\ref{t.main}. As we will see from the proof of (\ref{eqctrw}) and (\ref{eqctrwplus}), this is really an example where not leaving the set $\CC$ at all is ``responsible'' for almost 
all of the correlation between $\{\jumpch_0\in\CC\}$ and $\{\jumpch_t\in\CC\}$.

We first prove  (\ref{eqctrwmin}), which might already be in the literature somewhere, but we could not locate a reference. We start with the lower bound. In essence, the bound holds because $Y$ has probability of order $t^{(1-\beta)/2}$ to begin at a site of order at least $t^{1/2}$ from the vertex $1$. From such a site on the positive half-line, the walk has positive probability to remain positive for $t$ steps: indeed, at such sites the walk experiences an excess in leftward transition probability over rightward of order $t^{-1/2}$, so that, during $t$ steps, this imbalance provides a drift towards the origin totalling an order of $t^{1/2}$ steps. This drift is thus comparable to the Gaussian fluctuation of the particle during this period, and the particle remains to the right of the origin with positive probability. 
Rather than make this heuristic rigorous, we prove the lower bound in (\ref{eqctrwmin}) by invoking the Carne-Varopoulos bound (see \cite{carnevar,Var:speed} or \cite[Theorem 13.4]{LPbook}) which, in the case of  a reversible Markov chain $X$ having finite stationary measure  (and hence spectral radius $1$), asserts that
\begin{equation}\label{eqcarnevar}
  \Pb{ X(s) = y \md X(0)=x  } \leq 2 \sqrt{\frac{\pi(y)}{\pi(x)}} \exp \left\{ - \frac{d(x,y)^2}{2  s} \right\}
\end{equation}
for all $x,y \in S$, $s > 0$; here, $d(\cdot,\cdot)$ denotes graphical distance on $S$, and $\pi$ denotes the stationary measure. 
We apply this bound to the walk $Y$. Noting that $\tfrac{\pi(x)}{\pi(1)} \leq C x^{-\beta}$ for $x \in \N$, we find that,  if $s \in \{0,\ldots,t\}$,
and $x \in \N$ satisfies $t^{1/2} \sqrt{2} (\log t)^{1/2} \leq x \leq t^{1/2} \sqrt{2} (\log t)^{1/2}  + t^{1/2}$,
$$
  \Pb{ Y(s) = 1 \md Y(0)=x  } \leq C  \big( t^{1/2} \sqrt{2} (\log t)^{1/2}  + t^{1/2} \big)^{\beta/2} \exp \big\{ - 2 \log t \big\}. 
$$  
Sum this bound over $s \in \{0,\ldots,t\}$ to arrive at
$$
  \Pb{ \exists s \in \{0,\ldots,t\}: Y(s) = 1 \md Y(0)=x  } \leq C  t^{\beta/4 - 1} \big( \sqrt{2}  (\log t)^{1/2}  + 1 \big)^{\beta/2}. 
$$  
Recalling that $\beta < 2$, we find that, for $t$ high enough, the conditional probability, given that $Y(0)$ assumes any one of the $t^{1/2}$ values of $x$ described above, that $Y$ reaches $1$ before time $t$, is at most one-half. The probability that $Y(0)$ assumes some such value is at least 
$c  \big( t^{1/2} \sqrt{2} (\log t)^{1/2}  + t^{1/2} \big)^{-\beta} t^{1/2} 
= c t^{(1 - \beta)/2}  \big( \sqrt{2} (\log t)^{1/2} + 1 \big)^{-\beta}$. Thus, 
$$
  \Pb{ Y(s) \geq 1 \, \forall  s \in \{0,\ldots,t\} } \geq \tfrac{c}{2} t^{(1 - \beta)/2}  \big(  \sqrt{2} (\log t)^{1/2} + 1 \big)^{-\beta},
$$  
so that the lower bound in (\ref{eqctrwmin}) is verified.
 
We turn to the upper bound in (\ref{eqctrwmin}).
Let $Z$ denote the Markov chain on $\Z$ which shares its initial distribution with $Y$ and which evolves as a simple random walk. The two processes may be coupled so that, should $\jumpch(0)$ be positive, then  $\jumpch(m) \leq Z(m)$ for all $m$ at most the hitting time of $1$ by $\jumpch$; for this reason, it suffices  to establish the upper bound in (\ref{eqctrwmin}) for the process $Z$.

Let $a,t \in \N$. Let $Z_a:\{0,\ldots,t\} \longrightarrow \N$ denote simple random walk with $Z_a(0)=a$.
By \cite[Theorem 2.17]{LPW}, $\Pb{ \exists s \in \{ 0,\ldots,t\}: Z_a(s) = 0 } \leq 12 a t^{-1/2}$. Thus,
\begin{equation} \label{eqnupb}
\Pb{  Z(s) \geq 1 \, \, \forall \, \, 0 \leq s \leq t \md Z(0) = a  } \leq  12 a t^{-1/2} \, .
\end{equation}
 Multiplying the inequality resulting from (\ref{eqnupb}) by $\Pb{Z(0) = a}$, we sum over $a \in \N$ to obtain
\begin{equation}\label{eqzthalf}
\Pb{  Z(s) \geq 1 \, \, \forall \, \, 0 \leq s \leq t  } 
\leq 12 \sum_{a=1}^{t^{1/2}}  \Pb{Z(0)  = a} a t^{-1/2} \, + \, \Pb{Z(0) > t^{1/2}}\,.
\end{equation}
Note that there exists $C>0$ such that $\Pb{Z(0) = a} \leq C a^{-\beta}$ for $a \in \N$; thus, $\beta \in (1,2)$ implies that each of the two terms on the right-hand side of 
(\ref{eqzthalf}) is at most a constant multiple of $t^{(1 - \beta)/2}$. 
That is, the upper bound in (\ref{eqctrwmin}) holds for $Z$, as we sought to show. This completes the proof of (\ref{eqctrwmin}).

We now show (\ref{eqctrw}). Let $U$ be the first time that $\jumpch$ makes the jump $(-1,1)$, and $V$ be the first time that $\jumpch$ makes any of the three jumps $(-1,1)$, $(1,1)$, $(-1,-1)$. Obviously, $U \geq V$. The point of considering $V$ is that   $c_{-1,1}=c_{1,1}=c_{-1,-1}$ implies that $\Ps{\jumpch_V \in \CC \md \jumpch_{[0,V)}}=1/2$, and then the symmetry of the entire chain and the strong Markov property implies that $\Ps{\jumpch_t \in \CC \md V \leq t,\ \jumpch_{[0,V)}}=1/2$, as well.  On the other hand, $\big\{ \jumpch_s \in \CC \  \forall s \in [0,t]\big\}= \big\{\jumpch_0 \in \CC,\ U > t\big\}\supseteq  \big\{\jumpch_0 \in \CC,\ V > t\big\}$, and hence $\Pb{\jumpch_t \in \CC \md \jumpch_0 \in \CC,\ V > t}=1$. Therefore,
\begin{equation}\label{e.Y0Yt}
\begin{aligned}
\Pb{\jumpch_0,\jumpch_t \in \CC} &= \Pb{\jumpch_0\in \CC,\ V \leq t} \cdot \Pb{\jumpch_t\in \CC \md \jumpch_0
\in \CC,\ V \leq t}\\
&\quad + \Pb{\jumpch_0 \in \CC,\ V > t} \cdot \Pb{\jumpch_t \in \CC \md \jumpch_0 \in \CC,\ V > t} \\
& =  \big( \Pb{\jumpch_0\in \CC} - \Pb{\jumpch_0 \in \CC,\ V > t} \big) \cdot 1/2 + \Pb{\jumpch_0 \in \CC,\ V > t} \cdot 1\\
&\leq \frac{1}{2}\, \Pb{\jumpch_0\in \CC} +  \frac{1}{2}\, \Pb{\jumpch_0 \in \CC,\ U > t}\,.
\end{aligned}
\end{equation}
Using that $\Ps{\jumpch_t\in\CC}=\pi(\CC)=1/2$ for all $t\geq 0$, we get (\ref{eqctrw}).


By the second equality in (\ref{e.Y0Yt}), $\Pb{\jumpch_0,\jumpch_t \in \CC} = \tfrac{1}{4} + \tfrac{1}{2} \Pb{\jumpch_0 \in \CC,\ V > t}$.
The bound (\ref{eqctrwplus}) thus follows 
from the claim that for each $\epsilon > 0$, and for all $t > 0$ sufficiently high,
\begin{equation}\label{eqyzerocc}
\Pb{ \jumpch_0 \in \CC,\ V > t } \geq  \big( 1 - \epsilon \big) \Pb{ \jumpch_0 \in \CC,\ U > t }.
\end{equation}
The event on the right-hand side is simply $\{ \jumpch_s \geq 1 \  \forall s \in [0,t] \}$; the event on the left-hand side contains the event $\{ \jumpch_s \geq 2 \  \forall s \in [0,t]  \}$. Hence, it is enough to argue that the conditional probability of the latter event given the former tends to one in a limit of high $t$.  The event that $\jumpch_s \geq 1$ for all $s \in [0,t]$ and $\jumpch_s = 1$ for some such $s$ entails either that $Y$ remains positive for time $t/2$ after first reaching $1$ after time $0$, or that the same holds for the reversed chain, with time running backwards from $t$; either of these events has probability at most $C t^{-1/2}$ by (\ref{eqnupb}) applied for $a=1$. However, the event  $\jumpch_s \geq 1 \  \forall s \in [0,t]$ has probability at least $c t^{(1-\beta)/2 + o(1)}$ by the lower bound in (\ref{eqctrwmin}); this is much more probable under our hypothesis that $\beta < 2$, so that, given that $Z$ is strictly positive on $\{ 0,\ldots,t \}$, the conditional probability that $Z$ visits $1$ during this interval tends to zero in high $t$.  In this way, we obtain (\ref{eqyzerocc}) and thus (\ref{eqctrwplus}).

\subsection{Some cases of non-sharpness}

We first give an example where Theorem~\ref{t.main} is not at all sharp. Consider the process $\jumpch$ as above, and take $\CC$ to be the set of even positive integers. The conductances on the self-loops are set in such a way that $\Ps{Y_{t+1}\in 2\Z \md Y_t}=1/2$, regardless of $Y_t$.
Therefore, $\pi(2\Z)=1/2$ and $\pi(\CC)= 1/4$, and, using (\ref{eqctrw}) and (\ref{eqctrwplus}), then (\ref{eqctrwmin}), we find that, for $t \geq 1$,
\begin{align*}
\Pb{\jumpch_t\in\CC \md \jumpch_0\in \CC} &= \frac{1}{2}\, \Pb{\jumpch_{t-1} > 0 \md \jumpch_0 > 0}\\
&=  \frac{1}{2}\,\Big(1/2 + \Theta(1)\, \Pb{\jumpch_s > 0 \ \forall \, 0\leq s \leq t-1 \md \jumpch_0>0}\Big) \\
&= \frac{1}{4} + t^{\frac{1-\beta}{2} + o(1)}\,.
\end{align*}
Hence the correlation is polynomially large. On the other hand, $\Pb{\jumpch_s \in\CC \ \forall \, 0\leq s\leq t}$ is clearly exponentially small.

A more complicated but more natural example is given by Corollary~\ref{cor.LR} below.

\medskip

The previous subsection showed that Theorem~\ref{t.main} is sharp in the regime of polynomial decay. Examining the proof of the theorem, we have the feeling that this is not the case in the regime of superpolynomial decay. In particular, we have the following question.

\begin{question} 
Does the exponential pairwise decorrelation $d(t)=C\exp(-ct)$ in (\ref{e.decorr2}) for some event $\CC$ imply an exit time exponential decay $\Pb{\omega_s \in \CC \ \forall s\in [0,t]} < C'\exp(-c't)$, with $c',C'$ depending only on $c,C$ and $\Ps{\CC}=p$?
\end{question}

\ignore{
For reversible chains, it is well-known that the existence of some $c>0$ with 
\begin{equation}\label{e.onestep}
\Pb{\omega_0,\omega_1 \in \CC}-\Ps{\omega_0\in\CC}^2 \leq (1-c)\, \P[\CC](1-\P[\CC])
\end{equation}
for {\it all} static events $\CC$ is equivalent to having an absolute spectral gap for $T_1$ that depends only on $c$, hence it implies an exponential decorrelation in (\ref{e.decorr1}) and (\ref{e.decorr2}). As we discussed in the Introduction (and as shown, e.g., by the next subsection), in this case the answer to the above question is affirmative. However, if we have (\ref{e.onestep}) only for our single given $\CC$, that does not imply exponential pairwise decorrelation: one example is the process $Y$ on $\Z$ given above, with $\CC=\N^+$.\margin{-3}{G: I feel now that this paragraph is too obvious to include.}
}


\section{Applications to dynamical percolation}\label{s.dynperc}

{\bf Critical planar percolation} is a central object of probability theory and statistical mechanics; see \cite{Grimm,WWperc} for background. The best understood example is $\mathsf{Bernoulli}(1/2)$ site percolation on the triangular lattice, where the existence of a conformally invariant scaling limit is known. Roughly, if we consider percolation on the lattice of mesh $1/n$, and any collection $\Quad_1,\dots,\Quad_k$ of conformal images of rectangles, then the joint distribution of the left-right crossing events inside these $\Quad_i$'s has a limit that is conformally invariant. Moreover, one can define a continuum random limit object encoding all the macroscopic crossing events. See \cite{SchSmG} and the explanations and references there. In {\bf dynamical percolation}, every site is switching between being open and closed according to an independent exponential clock, in such a way that the stationary distribution on $\{0,1\}^{V_n}$ is critical percolation, where $V_n$ is the set of sites, and 0 represents ``closed'' and 1 represents ``open'' . This model has been studied from three closely related points of view: 
\begin{enumerate}
\hitem{i.1}{(1)} {\it How long does it take to change macroscopic crossings? Or, how {\bf noise sensitive} are the crossing events?} A reasonable guess is that this time-scale is given by the expected number of pivotal switches in the unit square (i.e., changes of the left-right crossing event) being of order one. 
Let $\Piv(n)$ be the expected number of sites in critical percolation in the unit square with mesh size $1/n$ that are pivotal for the left-right crossing; it is known for the triangular lattice that $\Piv(n) = n^{3/4 + o(1)}$ \cite{SW}. Then, in the stationary process, using Fubini's theorem and the linearity of expectation, the above time-scale is simply $n^2/\Piv(n) = n^{5/4 + o(1)}$. This guess, based merely on the expectation, has been confirmed by \cite{BKS,SS,GPS1}: if $t\,n^2/\Piv(n)=t\,n^{5/4+o(1)}$ sites are resampled in the unit square, then the correlation of crossing before and after the resampling is $t^{-2/3}$, up to constant factors, as $t\to\infty$ \cite[Eq.~(8.7)]{GPS1}. Similar, though slightly weaker, results have been proved for general conformal rectangles $\Quad$ with piecewise smooth boundary. Furthermore, even for critical (i.e., $\mathsf{Bernoulli}(1/2)$) bond percolation on $\Z^2$, where the existence of critical exponents such as that describing the growth of $\Piv(n)$ are not known, it follows from the proof of \cite[Corollary 1.2]{GPS1}, together with Eq.~(2.6) there, that, after resampling $t\,n^2/\Piv(n)$ edges, the correlation is at most $O(t^{-\alpha})$ for some $\alpha>0$. 
\hitem{i.2}{(2)} {\it On an infinite lattice, are there random times with exceptional behavior, e.g., with an infinite cluster? In other words, which events are {\bf dynamically sensitive}?} It was proved in \cite{SS} that there are exceptional times with an infinite cluster, and in \cite{GPS1} that their Hausdorff dimension is almost surely 31/36, and that such exceptional times also exist for critical dynamical bond percolation on $\Z^2$. A natural law on the infinite cluster that appears at exceptional times will be introduced and studied in \cite{HPS}. 
\hitem{i.3}{(3)} {\it In the unit square (or in another conformal rectangle), with mesh $1/n$ and a well-chosen rate $r(n)$ for the exponential clocks, is there a {\bf scaling limit of the process}, giving a Markov process on continuum configurations?} If we choose $r(n)=1/\Piv(n)=n^{-3/4+o(1)}$, then the expected number of pivotal switches
in the unit square during a unit time will be exactly 1, independently of $n$. 
It is proved in \cite{GPS2a, GPS2b} that, with this scaling, such a scaling limit does indeed exist on the triangular lattice. Additionally, it follows from the results of \cite{GPS1}, item~\iref{i.1} above, that the resulting Markov chain is ergodic;  in particular, the correlation decay for the unit square, in rescaled large time $t$, is $t^{-2/3}$, up to constant factors.
\end{enumerate}

We will have an application of Theorem~\ref{t.main} to the setup of items~\iref{i.1} and \iref{i.3}, and an application of Theorem~\ref{t.AKSz} to the setup of item~\iref{i.2}. Here is the first of these results:

\begin{corollary}\label{cor.LR} 
In dynamical critical site percolation on the triangular lattice or bond percolation on $\Z^2$, with mesh $1/n$ and rate $1/\Piv(n)$ for the clocks, consider the left-right crossing event $\CC$ in the unit square. There exist constants 
$\big\{ C_K: K \in \N \big\}$ such that, for each $K \in \N$ and for all $t>0$ and $n \in \N$,
$$
(1/4)^{\lceil 2t\rceil} \leq \Pb{ \omega_s \in \CC \, \forall \, s \in [0,t]   } \leq C_K t^{-K}.
$$
On the triangular lattice, it is known that $\Piv(n)=n^{3/4+o(1)}$, and the above bounds in $t$ also hold for the scaling limit of dynamical percolation.
\end{corollary}

Before starting the proof, let us emphasize that this corollary concerns a natural question that has exactly the kind of setup for which Theorem~\ref{t.main} is designed. Namely, in the finite $n$ version with the discrete-time chain (with sites being resampled one-by-one), the mixing time of the entire chain is $n^{2+o(1)}$ steps (it is just random walk on an $n^2$-dimensional hypercube), while the left-right crossing event $\CC$ decorrelates on the scale $n^{5/4+o(1)}$, as mentioned in item~\iref{i.1}. In the scaling limit of the chain, only the evolution of {\it macroscopic} crossing events is considered, so that $n^{5/4+o(1)}$ is the natural scaling factor needed to obtain this scaling limit. In particular, we are interested in the tail probability of exiting $\CC$ on this time scale, a question for which analysis based on the spectral gap of the entire chain would clearly be too crude. Moreover, it turns out that the limit chain does not have a spectral gap (something which is clear from the polynomial decorrelation $t^{-2/3}$), 
hence the classical exponential exit time results  \cite{AKSz} do not apply. One may nevertheless hope that at least
there would be a ``spectral gap restricted to $\CC$'', i.e., $\big( T_1(g\1_\CC), g\1_\CC \big) <  (1-c) (g,g)$ for some $c>0$, for all $g \in L^2(\{0,1\}^{V_n},\P_{1/2})$, which, similarly to the proof of Theorem~\ref{t.AKSz}, would imply an exponentially small upper bound in Corollary~\ref{cor.LR}. However, it is not hard to prove that $g=\1{\{\text{density of open bits is}\, > 1/2 + n^{-3/4+\eps}\}}$, with $\eps > 0$ fixed but small enough, is a counterexample. There are slightly more complicated counterexamples that make sense also in the scaling limit.

It is also interesting to note that Corollary~\ref{cor.LR}, despite being a consequence of Theorem~\ref{t.main}, provides a natural example in which the theorem by itself is not sharp: the correlation decay is polynomial, while the exit time tail is superpolynomial.

\proofof{Corollary~\ref{cor.LR}} We will work in the discrete lattice setting, i.e., with a fixed finite $n\in\Z^+$. All our results will hold uniformly in $n$, so that item~\iref{i.3} above implies that the results extend to the continuum scaling limit.

Firstly the upper bound. Let $L \in \N$, and decompose the unit square into $L$ vertical slabs with dimensions $1/L \times 1$ (the induced subgraphs in the slabs will not be exactly isomorphic to each other, but this is not a problem). For $s \geq 0$ and $i \in \{ 1,\ldots L \}$, 
let $A_i(s)$ denote the event that the $i^\textrm{th}$ such slab has an open left-right crossing at time $s$. We further write $\mathcal{A}_i(t)$ for the intersection of the events $A_i(s)$ over all $s \in [0,t]$. Clearly,
\begin{equation}\label{eqacap}
 \left\{ \omega_s \in \CC \, \forall \, s \in [0,t] \right\} \subseteq \bigcap_{i=1}^L \mathcal{A}_i. 
\end{equation}
We may apply Theorem~\ref{t.main} to bound $\Pb{\mathcal{A}_i}$. To do so, we need to have a correlation decay  
between the events $A_i(0)$ and $A_i(t)$. Indeed,
$$ 
\Pb{A_i(0),A_i(t)} -  \Pb{A_i(0)}^2  \leq  C'_L\, t^{-\alpha}
$$
holds for all $t\geq 0$ and some constant $C'_L$, with $\alpha=2/3$ in the case of the triangular lattice. This is simply the analogue for the slab of the decorrelation bound mentioned in item~\iref{i.1} above; if one does not want to optimize the constant $C'_L$, then the proof is identical to the one for the square case in \cite[Corollary 1.2]{GPS1}. Noting that the variance of the left-right crossing event in any one of the slabs is an $L$-dependent constant, Theorem~\ref{t.main} may be applied and yields that
\begin{equation}\label{eqtmaininf}
 \Pb{\mathcal{A}_i(t)} \leq t^{-\alpha + o(1)}, 
\end{equation}
for each $i \in \{1,\ldots,L\}$, where the $o(1)$ term depends on $L$. We now use (\ref{eqacap}), (\ref{eqtmaininf}) and the independence of dynamical percolation in disjoint regions to find that
$$
  \Pb{ \omega_s \in \CC \, \, \forall \, s \in [0,t] } \leq t^{-\alpha L + o(1)}, 
$$
and the upper bound follows.

For the lower bound, we need a basic tool that we call the  {\bf dynamical FKG inequality}. The next lemma is not the strongest possible form of such a result, but it will suffice for our purposes. Firstly, we need some notation.

A realization of dynamical percolation on a finite graph $G$ may be interpreted as a map 
$\omega: \R \longrightarrow \{ 0 , 1 \}^{V(G)}$. Let $\mathcal{S}$ denote the space of such maps; we will write $\omega_s(x)$ (for $s \in \R$ and $x \in V(G)$) for the value at time $s$ of $\omega$ in bit $x$.
 For static configurations, that is, elements $\eta \in \{ 0,1 \}^{V(G)}$, we consider the natural component-wise partial order: $\eta \preceq \eta'$ if{f} $\eta(x)\leq \eta'(x)$ for all $x\in V(G)$. We extend this to a partial order on $\mathcal{S}$ by writing $\omega \preceq \omega'$ if{f} $\omega_s \preceq \omega'_s$ for all $s \in \R$. A function $f: \mathcal{S} \longrightarrow \R$ is called increasing if $f(\omega) \leq f(\omega')$ whenever $\omega \preceq \omega'$. An event $C \subseteq \mathcal{S}$ is called increasing if $\One_C$ is increasing.
 
The standard Harris-FKG inequality for percolation (see \cite[Theorem 2.4]{Grimm}) says that increasing functions of static configurations are positively correlated. In particular, conditioning on an increasing static event makes the percolation configuration ``larger''. Our dynamical FKG inequality deals with conditioning on an increasing dynamical event:

\begin{lemma}\label{l.dynFKG}
Let $C \subseteq \mathcal{S}$ 
be an increasing event. Let 
$\P$ denote the law of dynamical percolation on 
$[0,\infty)$. Then there exists a coupling $\coupling$ 
of the laws $\P$ and $\Pb{ \cdot \md C }$ such that, denoting the two marginals by $\omega$
and $\omega^C$, we have that
$
 \coupling \big[   \omega_0 \preceq \omega^C_0   \big] = 1.
$
\end{lemma}

Assuming this lemma, if $\mathcal{A}(t)$ is the event that the square is crossed for all $s \in [0,t]$, then, for $k \in \N$,
\begin{equation}\label{eqank}
\Pb{\mathcal{A}((k+1)/2)} \geq \Pb{\mathcal{A}(k/2)} \Pb{\mathcal{A}(1/2)}.
\end{equation} 
Indeed, Lemma~\ref{l.dynFKG} implies that, conditionally on $\mathcal{A}(k/2)$, the distribution of the marginal of dynamical percolation at time $k/2$ stochastically dominates critical percolation. The event $\mathcal{A}(k/2)$ being conditionally independent of the subsequent evolution of dynamical percolation given the configuration at time $k/2$, we see that, conditionally on $\mathcal{A}(k/2)$, the distribution of dynamical percolation on $[k/2,(k+1)/2]$
stochastically dominates its unconditioned counterpart, whence (\ref{eqank}).

We claim now that $\Pb{\mathcal{A}(1/2)} \geq 1/4$. Indeed, the expected number of pivotal switches during a duration of one-half of scaled time is $1/2$, by Fubini's theorem, and, by symmetry, this remains the case conditionally on there being a left-right crossing at the start of this duration. So, by Markov's inequality, the conditional probability of having no pivotal switch during this time is at least $1/2$. The probability of a left-right crossing being $1/2$, we find that  $\Pb{\mathcal{A}(1/2)} \geq 1/4$. Thus, iterating~(\ref{eqank}) gives the lower bound in Corollary~\ref{cor.LR}. 
\qed

We still owe the proof of the dynamical FKG lemma that we used:

\proofof{Lemma~\ref{l.dynFKG}} Note that the law of $\omega$ may be constructed as follows. To each $x \in V(G)$, we associate a sequence $e_i(x)$, $i=1,2,\ldots$ of independent exponential mean one random variables, and an independent sequence $b_i(x)$, $i=1,2,\ldots$ of independent Bernoulli random variables. We set $\omega_0$ to be a uniform element in $\{0,1\}^{V(G)}$. We further set $\omega_t(x) = b_i(x)$ where $i \in \N^+$ is minimal subject to 
$\sum_{j=1}^i e_j(x) \leq t$. If no such $i$ exists, then $t < e_1(x)$; in this case, we set $\omega_t(x) = \omega_0(x)$. 

Write $\Omega^+$ denote the data $e_i(x)$ and $b_i(x)$ for $i \in \N^+$ and $x \in V(G)$. 
Note that $\Omega^+$ and $\omega_0$
comprise all of the data that specifies $\omega$. As such, we may denote an instance $\omega$ of dynamical percolation on $[0,\infty)$ in the form $(\omega_0,\omega^+) \in \{0,1\}^{V(G)} \times \Omega^+$. 

Suppose given an element $\omega^+ \in \Omega^+$.  Note that, if $\omega_0,\omega'_0 \in \{0,1\}^{V(G)}$
satisfy $\omega_0 \preceq \omega'_0$, then $(\omega_0,\omega^+) \in C$ implies that  $(\omega'_0,\omega^+) \in C$.
Hence, by the Harris-FKG inequality for a static configuration, the distribution of $\omega_0$ under $\Pb{ \cdot \md C, \omega^+ }$ 
stochastically dominates its distribution under $\Pb{ \cdot \md \omega^+ }$ (which is the uniform distribution).

Let $\mu_{C,+}$ denote the conditional distribution of $\omega^+$ given $C$. 
Then
\begin{equation}\label{eqaverage}
 \Pb{ \cdot \md C } = \int \Pb{ \cdot \md C, \omega^+ } \,  d\mu_{C,+}(\omega^+).
\end{equation}  
We have shown that, for all choices of $\omega^+$,  the conditional distribution of $\omega_0$ given $C$ and $\omega^+$ stochastically dominates the uniform distribution on $\{0,1\}^{V(G)}$. This statement remains true after the averaging in (\ref{eqaverage}).
Hence, we find that the law of $\omega_0$ given $C$ stochastically dominates its unconditioned law. 
\qed

\begin{question}\label{q.exp2/3}
What is the true decay of the probability for having a left-right crossing of the unit square  during $[0,t]$ in the scaling limit of dynamical percolation? We expect it to be $\exp(-t^{\beta+o(1)})$, with $\beta \in (0,1)$.
\end{question}

We now move to the application of Theorem~\ref{t.AKSz} to the study of exceptional times, which is item~\iref{i.2} in the list at the start of Section~\ref{s.dynperc}. Consider critical dynamical site percolation $(\omega_t)_{t\geq 0}$ on the infinite triangular lattice or bond percolation on $\Z^2$ (no scaling of space or time), and let 
$$
\EE:=\{t \in [0,\infty) : \omega_t\text{ has }0 \longleftrightarrow\infty\}
$$ 
be the set of exceptional times when the cluster of the origin is infinite. For any fixed time $t$, we have $\Ps{t\in\EE}=0$; hence $\EE$ has zero Lebesgue measure almost surely. However, as claimed in item~\iref{i.2}, it is almost surely nonempty, with Hausdorff dimension 31/36 in the case of the triangular lattice. A natural question is how long one has to wait to see the first exceptional time. It is answered by the following corollary which will also be an important tool for \cite{HPS} in studying the infinite clusters that appear in $\EE$.

\begin{corollary}\label{cor.radial} 
There exist $\infty > c_1 \geq c_2 > 0$ such that,
for critical dynamical site percolation on the infinite triangular lattice or bond percolation on $\Z^2$,
$$
\exp(-c_1\,t) \leq \Pb{\EE \cap [0,t] = \emptyset} \leq \exp(-c_2\,t)\,.
$$
\end{corollary}

\def\FET{\mathsf{FET}}
\def\fet{\mathsf{FET}}

In other words, the first exceptional time $\FET:=\min \EE$ has an exponential tail. Note here that, by \cite[Lemma 3.2]{HgPS}, the set $\EE$ is topologically closed, hence the minimum makes sense. Furthermore, $\Cl(\omega_0)$ is almost surely finite, hence it takes positive time until a bit in its boundary $\partial \Cl(\omega_0)$ first changes its status; thus $\FET>0$ almost surely. 

Naturally, the proof of the corollary will go through the finite approximations $\FET_R:=\inf\big\{t\in[0,\infty) : \omega_t\text{ has }0 \llra\p B_R(0)\big\}$. But first of all we need to prove that these times are actually approximations to $\FET$:

\begin{lemma}\label{l.fetR}
We have that $\fet_R\to\fet$ almost surely as $R\to\infty$.
\end{lemma}

For the proof, we will need another result proved in \cite[Lemma 3.2]{HgPS}:

\begin{lemma}\label{l.exctop} 
Almost surely, the set $\EE$ of exceptional times is disjoint from the set of times at which the status of a site is updated. 
\end{lemma}

\proofof{Lemma~\ref{l.fetR}} By item~\iref{i.2} above, $\mathcal{E} \cap (0,\infty) \not= \emptyset$, thus $\fet < \infty$, almost surely. The sequence $\fet_R$ is increasing and bounded above by $\fet$; thus, 
 there exists some random $\tau \in (0,\infty)$ such that $\fet_R\nearrow \tau$. Assume that $\tau < \fet$ with positive probability, which is to say that the cluster $\Cl_0(\omega_\tau)$ of the origin in $\omega_\tau$ is finite with positive probability. Almost surely, the set of flip times for any given site is a locally finite subset of $\R$, and the same holds if we take the union of the flip times over the finite set $S:=\Cl_0(\omega_\tau)\cup \partial \Cl_0(\omega_\tau)$. Therefore, on the event $|\Cl_0(\omega_\tau)|<\infty$, there exists some random $\eps>0$ such that the interval $(\tau-\eps,\tau+\eps)$ either has no flip times for $S$ (and hence the set $\Cl_0(\omega_t)$ remains unchanged), or it has a single flip time, $\tau$. On the other hand, we know that $\limsup_{s \nearrow \tau} \vert \Cl_0(\omega_s) \vert = \infty$ almost surely, which is consistent with the above only if $\tau$ is a flip time for a site in $\partial \Cl_0(\omega_\tau)$, closing at exactly time $\tau$, and if the reopening of this site creates a configuration in which $0 \llra \infty$. However, Lemma~\ref{l.exctop} shows that this circumstance has zero probability to occur at any time. We conclude that $\tau = \fet$ almost surely, which completes the proof. \qed

\proofof{Corollary \ref{cor.radial}}
Lemma~\ref{l.fetR} shows that the events $\{\EE_R \cap [0,t] =\emptyset\}$ increase to the event $\{\EE \cap [0,t] = \emptyset\}$. This means that bounds for $\EE_R$ that are uniform in $R$ imply the same bounds for $\EE$.

The lower bound is given simply by the probability that the bits (sites or bonds) neighbouring $0$ are closed during $[0,t]$.

For the upper bound, in order to apply Theorem~\ref{t.AKSz}, we need
the well-known fact that our Markov process has a spectral gap that is
uniform in $R$ (it is just continuous time random walk on the
hypercube $\{0,1\}^{B_R(0)}$, with unit rates on the edges), and also
need that $\Pb{\EE_R\cap [0,1]\not=\emptyset}>c>0$, uniformly in $R$.
For the case of the triangular lattice, this is part of \cite[Theorem
1.3]{SS}, while, for the case of $\Z^2$, part of \cite[Theorem
1.5]{GPS1}.
\qed

To conclude, let us point out the following interesting phenomenon. Although the results of \cite{GPS1} behave well under the $n$ vs.~$1/\Piv(n)=n^{-3/4+o(1)}$ space-time scaling, and hence it is not hard to show that, in the scaling limit of dynamical percolation in the full plane (mentioned in item~\iref{i.3} above), exceptional times when the ball of radius 1 is connected to infinity do exist and have Hausdorff dimension 31/36 a.s., the tail for the first exceptional time is expected to behave differently in the scaling limit than in the discrete case. If we try to think of the scaling limit process as unit-order regions flipping between being well-connected and not-at-all-connected (analogues of being open and closed in the discrete process) at a roughly unit rate, then it seems reasonable that, similarly to the discrete case, the tail behaviour of not having the connection from radius 1 to infinity is comparable to the obvious lower bound, the tail for not having a connection across the annulus between radii 1 and 2. However, we expect this annulus-crossing tail to be subexponential (see Question~\ref{q.exp2/3}), which would give a subexponential lower bound also here. The same issue from a different viewpoint is that we do not have any more the spectral gap that we needed in order to apply Theorem~\ref{t.AKSz}, hence there is no reason to hope for an exponential tail.

\ \\
{\bf Alan Hammond}\\
Department of Statistics, University of Oxford\\
\url{http://www.stats.ox.ac.uk/~hammond/}\\
\\
{\bf Elchanan Mossel}\\
Departments of Statistics and Computer Science, U.C. Berkeley, and\\ 
Faculty of Mathematics and Computer Science, Weizmann Institute of Science \\ 
\url{http://www.stat.berkeley.edu/~mossel/}\\
\\
{\bf G\'abor Pete}\\
Institute of Mathematics, Technical University of Budapest\\ 
\url{http://www.math.bme.hu/~gabor/}\\

\end{document}